\documentclass{article}
\usepackage{fullpage,amssymb}
\def\N{{{\rm I}\kern-.16em{\rm N}}}
\def\goback#1{\setbox0\hbox{#1}\kern-\wd0 \relax}
\def\ran{\mathop{\rm range}\nolimits}
\def\rank{\mathop{\rm rank}\nolimits}
\def\diag{\mathop{\rm diag}\nolimits}
\def\deg{\mathop{\rm deg}\nolimits}
\def\spa{\mathop{\rm span}\nolimits}
\def\ann{\mathop{\rm ann}\nolimits}
\def\C{{\rm C\kern-.48em\vrule width.06em height.57em depth-.02em
            \kern.48em}}
\def\nt{\noindent}
\begin{document}
\large

\title{Applications of the duality method to generalizations of
the Jordan canonical form}
\author{Olga Holtz  \\ 
\small Department of Mathematics \\
\small University of Wisconsin \\
\small Madison, Wisconsin 53706 U.S.A. \\
\small holtz@math.wisc.edu}
\date{}
\maketitle

We show how Ptak's duality method leads to short proofs of two 
extensions of the Jordan canonical form, viz.\ the normal form 
for a matrix over an arbitrary (not necessarily algebraically closed) 
field under similarity and the canonical form for a pair of matrices 
under contragredient equivalence. 

The duality method is summarized in the following. 

\nt {\bf Lemma.\/} {\em Let $V$ be a finite-dimensional space
over a field $F$, let $A : V\to V$ be a linear map, and $S\subset V$
be an $A$-invariant subspace of $V$. If $T\subset V^*$
is an $A^*$-invariant subspace of the dual $V^*$ of $V$ such that
\begin{eqnarray} s\in S, & \langle s, t \rangle =0  \quad \forall t\in T 
& \Longrightarrow \quad s=0, \label{prop1} \\
t\in T, & \langle s, t \rangle =0 \quad \forall s\in S & \Longrightarrow
\quad t=0, \label{prop2} \end{eqnarray}
then $V=S\dot{+}\ann(T)$ is an $A$-invariant direct sum decomposition 
of $V$, with $\ann(T)\colon= \{ v\in V : \langle v, t\rangle =0 
\quad \forall t\in T\}$ the annihilator of $T$. }

We give a proof for the sake of completeness.

{\bf{Proof.\/}} The condition~(\ref{prop1}) implies that
the sum $S+\ann(T)$ is direct. If $\dim T\geq \dim S$ and
$\{t_j\}_{j=1}^{\dim T}$ ($\{s_j\}_{j=1}^{\dim S}$) is a basis of $T$ ($S$),
then the matrix  $G\colon =(\langle s_i, t_j \rangle : i=1, \ldots, 
\dim S, \; j=1, \ldots, \dim T)$ has fewer rows than columns,
hence the equation $Gx=0$ has a nontrivial solution, so~(\ref{prop2}) 
fails. In other words, (\ref{prop2})~implies that 
$\dim T\leq \dim S$, hence $\dim \ann(T)\geq \dim V-\dim S$.
Thus, $V=S\dot{+}\ann(T)$. Since $T$ is $A^*$-invariant, $\ann(T)$
is $A$-invariant, which completes the proof. \hfill $\Box$ 

\section{The analogue of the Jordan form for an arbitrary field}

\nt {\bf Theorem 1.\/}  {\em Let $V$ be a finite-dimensional linear space
over a field $F$ and let $A : V\to V$ be a linear map. There 
exists a basis of $V$ such that the representation of $A$
with respect to that basis has the form
\begin{eqnarray}
& \diag(A_1, \ldots, A_p ), & \label{main} \\
& \mbox{where} & \nonumber \\
& A_i=\left(\begin{array}{ccccc} C_i & 0 & \cdots & 0 & 0 \\
 B_i & C_i & \cdots & 0 & 0 \\ \vdots & \vdots & \ddots & \vdots & \vdots \\ 
0 & 0 & \cdots & C_i & 0 \\ 0 & 0 & \cdots & B_i & C_i \end{array} \right),
\qquad  B_i=\left( \begin{array}{ccccc}0 & 0 & \cdots & 0 & 1 \\
0 & 0 & \cdots & 0 & 0 \\ \vdots & \vdots & \ddots & \vdots & \vdots \\
0 & 0 & \cdots & 0 & 0 \\ 0 & 0 & \cdots & 0 & 0 \end{array} \right)_{d_i
 \times d_i}\goback{$d_i\times$}, & \nonumber \\
& C_i=\left( \begin{array}{ccccccc} 0 & 0 & 0 & \cdots & 0 & 0 & a_{d_i} \\
1 & 0 & 0 & \cdots & 0 & 0 & a_{d_i-1} \\0 & 1 & 0 & \cdots & 0 & 0 & 
a_{d_i-2} \\ \vdots & \vdots & \vdots & \ddots & \vdots & \vdots & \vdots \\
0 & 0 & 0 & \cdots & 0 & 0 & a_{3} \\ 0 & 0 & 0 & \cdots & 1 & 0 & a_{2} \\
0 & 0 & 0 & \cdots & 0 & 1 & a_{1} \end{array} \right)_{d_i \times d_i} 
\goback{$d_i\times$}, &  \nonumber \\
& x^{d_i}-a_1x^{d_i-1}-\cdots-a_{d_i}\;\; {\hbox{is a prime 
in $F[x]$.}} & \nonumber 
\end{eqnarray}
This form is unique up to reordering of the blocks $A_1$, $\ldots$, $A_p$. }

{\bf{Proof.\/}} Since the space of all linear
maps on $V$ is finite-dimensional, there exists $k \in \N$ such that
$A^k\in \spa\{I,A,\ldots, A^{k-1}\}$, so $f(A)=\{0\}$, hence some monic
polynomial in $F[x]$ annihilates $A$. 

Let $f\in F[x]$ be the monic polynomial of minimal degree such that $f(A)=
\{0\}$ and let ${f=(f_1)^{k_1}\cdots (f_r)^{k_r}}$ be its decomposition into 
powers of distinct (monic) primes $f_i$, $i=1,\ldots, r$. 
Let $g_i\colon=\prod_{j=1, j\neq i }^r (f_i)^{k_i}$. Since $F[x]$
is a Euclidean domain and $\gcd(g_1, \ldots, g_r)=1$, it follows that
$g_1 h_1+\cdots +g_r h_r=1$ for some $h_1$, $\ldots$, $h_r\in F[x]$,
hence $v=h_1(A) g_1(A)v+\cdots+h_r(A) g_r(A)v$ for any $v\in V$. 
But $h_i(A) g_i(A)V \subseteq V_i\colon=\ker (f_i(A))^{k_i}$, so
$V=V_1+ \cdots + V_r$. Suppose $v\in V_i\cap \sum_{j\neq i} V_j$. 
As the polynomials $(f_i)^{k_i}$ and $\prod_{j\neq i}(f_j)^{k_j}$ are 
relatively prime, there exist $s_{i,1}$, $s_{i,2}\in F[x]$
such that $s_{i,1} (f_i)^{k_i}+s_{i,2}\prod_{j\neq i}(f_j)^{k_j}=1$, hence
$v=s_{i,1}(A)(f_i(A))^{k_i}v+s_{i,2}(A)\prod_{j\neq i}(f_j(A))^{k_j}v=0$, 
since $V_i=\ker (f_i(A))^{k_i}$ and $\sum_{j\neq i}V_j \subseteq \ker 
\prod_{j\neq i} (f_j(A))^{k_j}$. So, $V=V_1\dot{+}\cdots \dot{+}V_r$ is 
a(n $A$-invariant) direct sum decomposition of $V$. The arguments given so far
are standard.

Now show how to split the subspaces $V_i$. Let $\widetilde{V}$ stand for $V_1$, 
$\widetilde{A}$ for $A|_{V_1}$, $\widetilde{f}$ for $f_1$, $k$ for $k_1$,
 $d$ for $\deg f_1$. Since $f$ is the minimal polynomial annihilating
 $A$, $\widetilde{f}^{k}$ is the minimal polynomial annihilating 
$\widetilde{A}$, so there exists $v\in \widetilde{V}$ 
such that $w\colon =(\widetilde{f}(\widetilde{A}))^{k-1}\widetilde{A}^{d-1}
v\neq 0$. 

We claim that $w\notin \spa \{(\widetilde{f}(\widetilde{A}) )^{k-1}
\widetilde{A}^{j}v : j=0,\ldots, d-2 \}$. Indeed, if $w$ were in that span, 
it would imply $h(\widetilde{A}) (\widetilde{f}(\widetilde{A}))^{k-1}v=0$ 
for some polynomial $h$ of degree $d-1$. But any polynomial of degree $d-1$ 
is coprime to $f$, so there would exist a combination of $h$ and $f$ (with 
coefficients from $F[x]$) equal to $1$, which would yield $(f(\widetilde{A}
))^{k-1}v=0$, contradicting $w\neq 0$. Hence the claim follows.

So, there exists $v'\in \widetilde{V}^*$ such that
$$ \langle (f(\widetilde{A}))^{k-1} \widetilde{A}^j v,v'\rangle\cases{ 
=0 & if $j=0, \ldots, d-2$ \cr \neq 0 & if $j=d-1$.} $$ 
Let \begin{eqnarray*}
&& W_1\colon = \spa\{ (f(\widetilde{A}))^{i_1-1} \widetilde{A}^{i_2-1}v: 
i_1=1,\ldots, k, \;\; i_2=1, \ldots, d\}, \\
&& W_1'\colon = \spa\{ (f(\widetilde{A}^*))^{i_1-1} (\widetilde{A}^*)^{i_2-1}
v' : i_1=1,\ldots, k, \;\; i_2=1, \ldots, d\}. 
\end{eqnarray*}
Notice that 
$$ g_{(i_1,i_2),(j_1,j_2)}\colon= \langle (f(\widetilde{A}))^{i_1-1} 
\widetilde{A}^{d-i_2}v, (f(\widetilde{A}^*))^{k-j_1} (\widetilde{A}^*)^{j_2-1}
v' \rangle \neq 0$$ only if $(i_1,i_2) \preceq 
(j_1,j_2)$ (in lexicographic order). So, the $kd\times kd$-matrix
$(g_{(i_1,i_2),(j_1,j_2)} : i_1, j_1=1, \ldots, k, \; i_2, j_2=1, \ldots, d)$
is upper triangular with nonzero diagonal elements, hence, by the Lemma,
$\widetilde{V}=\dot{+}\ann(W'_1)$ is an $\widetilde{A}$-invariant direct 
sum decomposition of $\widetilde{V}$.
The matrix representation of $\widetilde{A}|_{W_1}$ with respect to the basis
$( (f(\widetilde{A}))^{i_1-1} \widetilde{A}^{i_2-1}v : i_1=1,\ldots, k, 
i_2=1, \ldots, d) $
ordered lexicographically is one of the diagonal blocks in~(\ref{main})
with $d_i=d$ and $\widetilde{f}(x)=x^d-a_1x^{d-1}-\cdots-a_d$.
 
Splitting the spaces $\ann(W'_1)$, $V_2$, $\ldots$, $V_r$ in the same way 
as above, we obtain a direct sum $V=W_1\dot{+}\cdots \dot{+}W_p$ of 
$A$-invariant 
indecomposable subspaces and a basis in each so that the matrix representation
of $A$ with respect to the concatenation of the bases of $W_i$'s has 
the form~(\ref{main}).    
       
Since the minimal polynomial $f$ of $A$ is unique, the (monic)
prime factors $f_i$ and the powers $k_i$ with which they occur in $f$
are determined uniquely. Let 
$$ n^i_j\colon =\dim \ker (f_i(A))^j=\sum_{W_l\subseteq \ker (f_i(A))^{k_i}}
\min (\dim W_l, j\,\deg f_i), \qquad i=1, \ldots, r, \; \;
j=1, \ldots, k_i. $$
Then $\Delta n_j^i\colon = n_{j+1}^i-n_j^i$ is the number of blocks for $f_i$ 
of order greater than $j\cdot \deg f_i$, so the number 
of blocks of order $j\cdot \deg f_i$ equals $-\Delta^2 n_{j-1}^i/\deg f_i=
(\Delta n_{j-1}^i- \Delta n_j^i)/\deg f_i$. Since the numbers $n_j^i$ are 
uniquely determined by the map $A$, this completes the proof of the uniqueness 
of~(\ref{main}). \hfill $\Box$

{\bf{Remarks.\/}} 1. The arguments in the two preceding paragraphs 
are variations of those due to de~Boor~\cite{dB}. 2. If $F$ is 
algebraically closed, the polynomials $f_i$ are of degree $1$, 
so~(\ref{main}) becomes the Jordan normal form of $A$. 3. In the proof above, 
all the factors of the minimal polynomial are treated in the same way in 
contrast to the
proof in~\cite{Pt} where the canonical splitting is first given for the
nilpotent part of $A$ and then follows for all other parts by shifting
$A$ by an eigenvalue $\lambda$~(for that completion of the proof in~\cite{Pt}, 
see~\cite{dB}). 4. Theorem~1 is classical and can be found,
e.g., in~\cite[pp.~92--97]{J}. In the sequel, we refer to a matrix in 
the form~(\ref{main}) as being in the {\em Jordan normal form for the 
field $F$\/}, and as the Jordan normal form of the operator $A$.
    
\section{The canonical form under contragredient equivalence}

Two pairs of matrices, $(A,B)$ and $(C,D)$, are called contragrediently 
equivalent if $A,C\in F^{m\times n}$, $B,D\in  F^{n\times m}$, and 
$A=SCT^{-1}$, $B=TDS^{-1}$  for some invertible $S\in F^{m\times m}$, 
$T\in F^{n\times n}$.

The problem of classification of pairs of matrices under contragredient 
equivalence can be restated as follows. Given an $n$-dimensional linear space
 $V$ and an $m$-dimensional linear space $W$ and linear
maps $A : V\to W$, $B : W \to V$, choose bases of $V$ and $W$ so that
the pair $(A,B)$ has a simple representation with respect to these bases.

\nt {\bf Theorem 2.\/} {\em Let  $V$, $W$ be finite-dimensional linear spaces 
over a field $F$ and 
let  $A : V\to W$, $B : W \to V$ be linear maps. There exist bases of $V$ and 
$W$ such that, with respect to those bases, the pair $(A,B)$ has the 
representation
\begin{equation}
\left( \diag(I, A_1, \ldots, A_p,0 ), \quad 
\diag (J_{AB}, B_1, \ldots, B_p, 0) \right) \label{*}
\end{equation}
where $J_{AB}$ is the nonsingular part of the Jordan form of $AB$,
$A_i, B_i\in F^{m_i\times n_i}$, $|m_i-n_i|\leq 1$, and 
$$(A_i,B_i)\in \{ (\left( \begin{array}{cc} I_{m_i-1}&  0 \end{array}\right), 
\left( \begin{array}{c} 0 \\ I_{m_i-1}\end{array} \right) ), \;
(\left( \begin{array}{c} 0 \\ I_{m_i-1}\end{array} \right), 
\left(\begin{array}{cc} I_{m_i-1}& 0 \end{array} \right)),\;
(I_{m_i},J_{m_i}), \; (J_{m_i},I_{m_i}) \}$$
where $J_k$ denotes the $k\times k$-matrix with ones on the first 
subdiagonal and zeros elsewhere.
The representation~(\ref{*}) is unique up to the order of the pairs of blocks
$(A_i,B_i)$, $i=1,\ldots, p$. Two pairs $(A,B)$ and $(C,D)$ are 
contragrediently equivalent if and only if $AB$ is similar to $CD$ and
\begin{eqnarray}
& \rank A=\rank C, \; \rank BA=\rank DC, \; \ldots,
 \; \rank (BA)^t=\rank (DC)^t, &\nonumber \\ 
& \rank B=\rank D, \; \rank AB=\rank CD, \; \ldots,
\; \rank (AB)^t=\rank (CD)^t, & \label{seq} \\ 
& t\colon=\min \{m,n\}.& \nonumber 
\end{eqnarray}   }
  
{\bf{Proof.\quad Step 1.\/}} By~Theorem~1 of~\cite{Pt} (whose
proof holds over an arbitrary field), there exist $V_1$ ($W_1$) and
$V_2$ ($W_2$) such that $BA$  ($AB$) is invertible on $V_1$ ($W_1$)
and nilpotent on $V_2$ ($W_2$) and $V=V_1\dot{+}V_2$ 
($W=W_1\dot{+}W_2$). Moreover, $V_1=\ran (BA)^r$ ,$V_2 = \ker (BA)^r$, 
$W_1= (AB)^r$, $W_2 = \ker (AB)^r$ for some $r\in \N$. 
If $x\in V_1$, then $x=(BA)^r y$ for some $y\in V$, hence $(AB)^r Ay=Ax$, 
that is, $Ax\in W_1$. Analogously, $By\in V_1$ whenever $y\in W_1$.
So, $V=V_1 \dot{+} V_2$, $W=W_1\dot{+} W_2$,
$A$ maps $V_i$ to $W_i$, $B$ maps $W_i$ to $V_i$ for $i=1,2$.

If $x\in V_2$, then $(AB)^r Ax=0$, so $Ax\in W_2$.
If $x\in V_1$ and $Ax=0$, then $BAx=0$, therefore, $x=0$, since $BA$
is invertible on $V_1$. So, $A$ induces a one-one map from $V_1$ to $W_1$.
Likewise, $B$ induces a one-one map from $W_1$ to $V_1$. So, $V_1$ 
and $W_1$ have the same dimension and the induced maps are also onto.

This step of the proof not only uses Theorem~1 of~\cite{Pt}, but
also parallels it.

 Now one can choose bases of $V_1$ and $W_1$ so that 
$A|_{V_1}$ is the identity
matrix and $B|_{W_1}$ is in Jordan normal form (which is the nonsingular
part of the Jordan normal form of $AB$). 

{\bf{Step 2.\/}} The spaces $V_2$ and $W_2$ are further split as follows.
 Let $l$ be the length 
of the longest nonzero product of the form $\cdots ABA$ or $\cdots BAB$.
Call such a product $C$ and suppose it ends in $A$. Pick $x\in V_2$ so that
$Cx\neq 0$ and form the sequence $x$, $Ax$, $BAx$, $\ldots$, $Cx$, whose
elements are alternately in $V_2$ and $W_2$. Let $V_3$ ($W_3$) be the span 
of the elements of the sequence belonging to $V_2$ ($W_2$).  

If $l$ is even, then $\dim V_3=\dim W_3+1=1+l/2$. Pick $x'\in V_2^*$ so that
$\langle Cx,x'\rangle\neq 0$. Form the sequence 
$x'$, $B^*x'$, $\ldots$, $A^*B^*x'$, $\ldots$, $C^* x'$. 
Let $V_4$ ($W_4$) be the annihilator in $V_2$ (in $W_2$) of the elements 
of the sequence that lie in $V_2^*$ ($W_2^*$). The $(1+l/2)\times
(1+l/2)$-matrix $(\langle (BA)^{i-1}x,(A^*B^*)^{1+l/2-j}x'\rangle : i,j=1,
\ldots, 1+l/2)$ is upper triangular with nonzero diagonal entries, hence,
by the Lemma,  $V_2=V_3\dot{+}V_4$. This argument is exactly the same as 
the corresponding argument in~\cite{dB}.

Analogously,  $W_2=W_3\dot{+} W_4$. Moreover,  $A$ 
maps $V_i$ to $W_i$, $B$ maps $W_i$ to $V_i$, $i=3,4$, and the pair 
$(A|_{V_3},B|_{W_3})$ has the form 
$$ ( \left(\begin{array}{cc} I_{l/2} & 0 \end{array}\right),
\left( \begin{array}{c} 0 \\ I_{l/2} \end{array} \right) ).$$   

If $l$ is odd, then $\dim V_3=\dim W_3=(1+l)/2$, and the above construction
gives $V_2=V_3\dot{+}V_4$,  $W_2=W_3\dot{+} W_4$ with $A$ mapping $V_i$ 
to $W_i$, $B$ mapping $W_i$ to $V_i$, $i=3,4$, the pair $(A|_{V_3}, B|_{W_3})$
having the form $(I_{(1+l)/2},J_{(1+l)/2})$.

If $C$ ends in $B$, then $(A|_{V_3}, B|_{W_3})$ has the form 
$$ (\left( \begin{array}{c} 0 \\ I_{l/2} \end{array} \right) ,\left(
\begin{array}{cc} I_{l/2} & 0 \end{array}\right) )\qquad \mbox{or} \qquad
(J_{(1+l)/2},I_{(1+l)/2}). $$

This step of the proof parallels, with necessary modifications, Theorem~2 
of~\cite{Pt}.

The problem is now reduced to splitting $V_4$ and $W_4$ in the same way.
The splitting process ends at the $j$-th stage if $A|_{V_{2j}}=0$
and $B|_{W_{2j}}=0$.

Thus one obtains the canonical form~(\ref{*}). It is completely determined by
the nonsingular part of the Jordan form of $AB$ and the ranks $\rank (A)$, 
$\rank (BA)$, $\rank (ABA)$, $\ldots$, $\rank (B)$, $\rank (AB)$, $\rank 
(BAB)$, $\ldots$. Since the rank of any such product equals the size of
$J_{AB}$ if the length of the product exceeds $2\min\{m,n\}$, the
infinite sequences above can be terminated at $(BA)^{\min\{m,n\}}$,
$(AB)^{\min\{m,n\}}$. It follows that 1) the representation~(\ref{*}) is 
unique up to the order of the pairs of blocks and that 2) two pairs $(A,B)$ 
and $(C,D)$ are contragrediently equivalent if and only if $AB$ is similar 
to $CD$ and~(\ref{seq}) holds. \hfill $\Box$

{\bf{Remarks.\/}} 1. Ptak's duality method was rediscovered by 
I.~Kaplansky~\cite{K}, who also described how to derive the canonical 
form~(\ref{*}). The same form was first published by N.~T.~Dobrovol'skaya
 and V.~A.~Ponomarev~\cite{DP}.  
J.~Gelonch and P.~Rubi\'{o}~i~Diaz~\cite[Theorem~2]{GR} proved that
the pair $(A, B)$ can be represented as 
$$\left( \diag (A_1, \ldots, A_q), \quad  \diag (B_1, \ldots, B_q) 
\right)$$
where $A_i$ and $B_i^*$ are of the same size and 
$$ (\dim\ker A_i,\dim\ker B_i)\in \{(0,1), (1,0) \}\quad 
\mbox{unless} \quad A_i=0, \; B_i=0.$$ R.~Horn and D.~Merino derived the 
canonical form~(\ref{*}) in~\cite[Theorem~5]{HM}. All the derivations 
(in~\cite{DP}, \cite{K}, \cite{GR}, and~\cite{HM}) were for the field $\C$. 2. 
Observe that the canonical form of the pair $(I,A)$ under
contragredient equivalence is $(I, J_A)$, where $J_A$ is the 
Jordan normal form of $A$. This and many other applications of the 
canonical form~(\ref{*}) are discussed in~\cite{HM}.


\begin{thebibliography}{11}

\bibitem{dB} C. de Boor, \smallskip On Ptak's derivation of the Jordan
canonical form, {\em Linear Algebra Appl.\/}
this issue.

\bibitem{DP} N. T. Dobrovol'skaya and V. A. Ponomarev, \smallskip A pair
of counter operators, {\em Uspehi Mat.\ Nauk\/} 20 : 80--86 (1965).  

\bibitem{GR} J. Gelonch and P. Rubi\'{o} i Diaz, \smallskip
Doubly multipliable matrices, {\em Rend.\ Istit.\ Mat.\ Univ.\
Trieste\/} 24 (vol.1--2) : 103--126 (1992).

\bibitem{HM} R. Horn and D. Merino, \smallskip Contragredient equivalence:
a canonical form and some applications, {\em Linear Algebra Appl.\/} 
214 : 43--92 (1995).

\bibitem{J} N. Jacobson, \smallskip Lectures in abstract algebra, 
New York, Springer-Verlag, 1975.

\bibitem{K} I. Kaplansky, \smallskip Private communication.

\bibitem{Pt} V. Ptak, \smallskip A remark on the Jordan normal form of matrices,
{\em Linear Algebra Appl.\/} this issue.

\end{thebibliography}
\end{document}